\documentclass[centertags,leqno,11pt]{article}

\usepackage{amsmath}
\usepackage{amssymb}
\usepackage{latexsym}
\usepackage{amsxtra}
\usepackage{amscd}
\usepackage{theorem}

\usepackage{graphics}
\usepackage{epic}

\theoremstyle{change}

{\theorembodyfont{\slshape}

\newtheorem{thm}{Theorem.}[section]
\newtheorem{cor}[thm]{Corollary.}
\newtheorem{lem}[thm]{Lemma.}
\newtheorem{prop}[thm]{Proposition.}

\newtheorem{defn}[thm]{Definition.}}

{\theorembodyfont{\rmfamily}

}

\newcommand{\proof}{\noindent {\bf Proof:\ }}
\newcommand{\Endproof}{\hspace*{\fill} $\Box$ \vspace{1ex} \noindent }

\makeatletter
\renewcommand{\subsection}{\@startsection{subsection}{2}%
{\z@}{-3.25ex plus -1ex minus-.2ex}{-1em}{\bf}} \makeatother

\newcommand{\PP}{\mathbb{P}}
\newcommand{\ZZ}{\mathbb{Z}}
\newcommand{\CC}{\mathbb{C}}

\newcommand{\QQ}{\mathbb{Q}}
\newcommand{\NN}{\mathbb{N}}
\newcommand{\FF}{\mathbb{F}}
\renewcommand{\AA}{\mathbb{A}}
\newcommand{\GG}{\mathbb{G}}

\newcommand{\V}{\mathcal{V}}

\newcommand{\bQl}{{\bar{\QQ}_\ell}}

\newcommand{\tame}{{\rm tame}}

\newcommand{\Hyp}{{\rm Hyp}}

\newcommand{\Unip}{{\rm {\bf U}}}

\newcommand{\p}{\mathfrak{p}}
\newcommand{\GL}{{\rm GL}}

\newcommand{\Gal}{{\rm Gal}}

\newcommand{\Spec}{{\rm Spec\,}}

\newcommand{\an}{{\rm an}}
\newcommand{\red}{^{\rm red }}

\renewcommand{\red}{{\rm red}}

\newcommand{\To}{\;\longrightarrow\;}

\newcommand{\Frob}{{\rm Frob}}

\newcommand{\G}{{\rm  G}}

\newcommand{\SO}{{\rm SO}}

\numberwithin{thm}{section}
\theoremstyle{plain}

\begin{document}

\title{Rigid local systems and potential automorphy: The $G_2$-case. }

\author{Michael Dettweiler\footnote{The author gratefully acknowledges financial support 
from the DFG Heisenberg Grant DE-1442}}
\maketitle

%\begin{document}

\begin{abstract} 
For $s\in \PP^1(\QQ)\setminus \{0,1,\infty\}, $ we study the family of compatible systems of  
Galois representations $$(\rho_\ell^s(3):\Gal(\bar{\QQ}/\QQ)\to
G_2(\QQ_\ell))_{\ell}$$  introduced 
in \cite{dk}, where $G_2(\QQ_\ell)\leq \GL_7(\QQ_\ell)$ is the simple group of type $G_2.$
We 
prove that, under some mild condition on $s,$ the image of the Tate twisted Galois representation $\rho_\ell^s(3)$ coincides with 
$G_2(\ZZ_\ell)$ for all $\ell$ up to a set of primes having density zero, and, using 
a potential automorphy criterion of Barnet-Lamb, Gee, Geraghty and Taylor
\cite{BLGGT}, we show 
that   $\rho_\ell^s(3)$ is potentially automorphic for all $\ell.$ \end{abstract}
%\tableofcontents 

%------------------------------------------------------

%\addcontentsline{toc}{section}{Introduction}
%\setcounter{page}{1}
\section*{Introduction}\label{Introduction}

Rigid local systems have a rich history, starting with Riemann's study of
Gauss'
hypergeometric differential  equations. N. Katz gave a 
unifying approach to all rigid local systems, using the theory 
of the middle convolution \cite{Katz96}. 
 In \cite{DR07} and \cite{dk}, we use Katz' theory of the middle convolution
in order to prove the existence of a family of motives having generically 
a motivic Galois group $G_2.$
The $\ell$-adic realization of this family is  an $\ell$-adic lisse  sheaf $\V_\ell$ of rank $7$  on 
$T=\PP^1_{\ZZ[\frac{1}{2\ell}]} \setminus \{0,1,\infty\}$ of weight $6.$ 
The lisse sheaf $\V_\ell$ corresponds to 
a continuous representation 
$$ \rho_\ell: \pi_1(T,\bar{s})\To \GL_7(\QQ_\ell),$$ where
$\bar{s}$ is a geometric point. If $F$ is a field and if $s\in T(F),$ one may consider 
the {\it specialization of $ \rho_{\ell}$ at $s,$} defined as the composition 
$$ \rho_\ell^s:  \Gal(\bar{F}/F)\tilde{\To} \pi_1(s,\bar{s})\To \pi_1(T,\bar{s})\stackrel{\rho_{\V_\ell}}{\To}
\GL_7(\QQ_\ell),$$ where $\bar{s}$ is a geometric point extending $s,$ and where 
the middle map is given by functoriality of the functor $\pi_1.$    It is shown \cite{dk} that if $s$ is 
a rational point of $T$ for which there exist two odd primes $p,q\not=\ell$ such that 
$s$ is $p$-adically close to $\infty$ and $q$-adically close to $1,$ then the specialization 
$ \rho^s_\ell(3): G_\QQ\to G_2(\QQ_\ell)\leq \GL_7(\QQ_\ell)$ has an image which is Zariski dense in 
the group $G_2(\QQ_\ell)$ (here we use  a third Tate twist  in order to 
obtain the monodromy group as a subgroup of the orthogonal group).   
It is natural to ask, whether these Galois representations correspond under 
the Langlands correspondence to an automorphic representation of $\GL_7(\AA).$  
Our main result is the following  (see Thm \ref{thmhaupt}):

\begin{thm} \label{thm01} Let  $s\in T(\QQ)$ is such that there exist two 
different primes $p,q\not= 2$ such that $\nu_p(s)>0$ and $\nu_q(1-s)<0.$
Then  the  Galois 
representations $\rho_\ell^s(3)$ are potentially automorphic for all 
$\ell.$ 
\end{thm}

The proof of this result relies on the compatibility and irreducibility of the family of 
specialized Galois representations $(\rho_\ell^s(3))$ (Thm.~\ref{dkkor}) and on 
the Hodge-Tate-regularity of the underlying motives  (Lemma~\ref{lemdk}).
Using self-duality, we can then apply the recent potential modularity result
\cite{BLGGT}, Thm.~A. 

The author thanks Nick Katz  for mentioning the possible application of Larsen's theorem to the problem of 
generic largeness of the Galois images in the $G_2$-case.

\section{Notation and definitions} \label{secdef}

Let us first set up some notation and definitions: 
If $k$ is a field, then we write $\bar{k}$ for an algebraic closure of $k$ and we set
$\Gamma_k=\Gal(k^{\rm sep}/k).$ If $F$ is a number field and $\omega$ is a finite prime of $F,$
then the completion of $F$ with respect to $\omega$ is denoted by $F_\omega,$ 
and the inertia subgroup of $\Gamma_{F_\omega}$ is denoted by $I_\omega.$  
The tame inertia group at $\omega$  is denoted by $I_\omega^\tame$
(quotient of $I_\omega $ by its $\ell$-Sylow subgroup,
where $\ell$ is the characteristic of $\omega$). 
Remember  there is a character 
$\omega_\ell: I_\omega^\tame \to \FF_\ell^\times,$ obtained by adjoining 
 the $\ell-1$-th roots of unity of $\ell.$  Usually, a primitive $k$-th root of unity 
 in some field is denoted by $\zeta_k.$ 
 
 If $\V$ is a lisse sheaf on a connected scheme $X$ and if $\bar{s}$ is a geometric point of 
 $X,$ then the corresponding 
 monodromy representation is denoted by $\rho_\V:\pi_1(X,\bar{s})\to \GL(\V_{\bar{s}}).$

\begin{defn} Let $R$ be a subring of $\CC,$ let $t$ be the standard parameter of $\AA^1_R\subseteq \PP^1_R$ and identify
$\CC((\frac{1}{t}))$ with ${\mathcal O}_{\PP^1_\CC,\infty(t)}.$ Then one may
identify
$$ \pi_1(\Spec \CC((\frac{1}{t})))\simeq \lim_{\stackrel{\leftarrow}{N}}
\Gal(\CC((\frac{1}{T^{\frac{1}{N}}}))/\CC((\frac{1}{t})))
\simeq \prod_p \ZZ_p.$$
Then the restriction of the monodromy representation $\rho_\V$ of a lisse sheaf $\V$ 
on $\AA^1_R\setminus \{0,1\}=
\PP^1_R\setminus \{0,1,\infty\}$ to the spectrum of $\CC((\frac{1}{t}))$ 
is called {\it the local monodromy of} $\V$ at $\infty.$ This notion
extends in the obvious way to the notion of local monodromy at $0,1.$ 
We call the local monodromy {\it of type}
$$ \Unip(n_1)\oplus\cdots\oplus  \Unip(n_k)$$
if it decomposes into $k$ indecomposable unipotent representations
of lengths $n_1,\ldots,n_k$ (resp.). 
\end{defn}

\section{A lisse sheaf of type $G_2$ and a family of $G_2$-motives}\label{secrevv}
Let $\ell$ be an odd prime, let 
 $R=\ZZ[\frac{1}{2\ell}],$ and let $T:=\AA^1_R\setminus\{0,1\}=\Spec(R[x][\frac{1}{x(x-1)}]).$
 The equation
\begin{equation}\label{eq-1} Y^2=\prod_{i=1}^7(X_{i+1}-X_i)\prod_{i=1,3,5,7}X_i\prod_{i=1,2,4,6}(X_i-1)
\end{equation} defines a smooth subscheme $\Hyp$ of $\GG_{m,R}\times (\AA^7_{X_1,\ldots,X_7,R}\setminus D),$
where $D$ is the vanishing locus of the right hand side of Equation~\eqref{eq-1}. 
Let $\sigma$ denote the involutory automorphism of $\Hyp,$ defined
by sending $Y$ to $-Y,$ and let $\p$ be the formal linear
combination $\frac{1}{2}(1+\sigma),$ viewed as an element in the group
ring of $\langle \sigma \rangle.$ 
Let $\pi:\Hyp \to T=\AA^1_{X_7,R}\setminus \{0,1\}$ denote the composition of the 
projection of $\Hyp$ onto $\AA^7_R$ followed by the projection 
onto the $7$-th coordinate. 
 By Deligne's work on  the Weil
conjectures (Weil II, \cite{DeligneWeil2}), the higher direct image with compact supports
$R^6\pi_!(\QQ_\ell)$ is mixed of weights $\leq 6.$ 
Moreover, the element $\p$ operates idempotently on $R^6\pi_!(\QQ_\ell)$
 and therefore cuts out a subsheaf from $R^6\pi_!(\QQ_\ell),$
denoted by 
$\p(R^6\pi_!(\QQ_\ell)).$ As shown in \cite{DR07}, it is a consequence of 
Katz' work on rigid local systems (\cite{Katz96}, Chap.~8), that the 
 weight-$6$-quotient 
$$\V_\ell:=W^6\left(\p(R^6\pi_!(\QQ_\ell))\right)$$ of $\p(R^6\pi_!(\QQ_\ell))$
is lisse on $T.$ Note that since a Tate twist diminishes the weight by $2,$ the Tate twisted sheaf 
$\V_\ell(3)=W^6\left(\p(R^6\pi_!(\QQ_\ell))\right)(3)$ has weight zero.  We set
$ \rho_\ell:=\rho_{\V_\ell}.$
The following motivic interpretation for $\rho_\ell$ is proved
in \cite{DR07}, Cor. 2.4.2:

\begin{lem} \label{cormotg4} The following holds:
\begin{enumerate}
\item There 
exists a smooth and projective scheme $X$  over 
$\AA^1_{\QQ}\setminus \{0,1\}$ and  an open embedding of 
$j:\Hyp_\QQ \to X$
 such that 
$$D=X\setminus \Hyp_\QQ=\bigcup_{i\in I}D_i$$ 
is a strict
normal crossings divisor over $\AA^1\setminus \{0,1\}.$ 
The involutory 
automorphism $\sigma$ of $\Hyp$ (given by $Y\mapsto -Y$)
 extends to an automorphism $\sigma$ of 
$X.$
 \item Let $\coprod_{i\in I}D_i$ denote the disjoint union 
of the components of $D$ and let 
$$\pi_X:X\to \AA^1_{\QQ}\setminus \{0,1\}\quad{\rm and}\quad 
 \pi_{\coprod D_i}: \coprod_{i\in I}D_i\to \AA^1_{\QQ}\setminus \{0,1\}$$ 
denote the structural 
morphisms. 
Then 
$$ \V_\ell|_{\AA^1_\QQ\setminus \{0,1\}}\simeq  \p\left[\ker\left(R^6(\pi_{X})_*(\bQl) \to R^6(\pi_{\coprod D_i})_*(\bQl)\right)\right],$$
where $\p$ denotes the formal sum $\frac{1}{2}(\sigma - 1).$ 
\end{enumerate}
\end{lem}

The following result is proved in  \cite{DR07}, Theorem~2.4.1:

\begin{prop}\label{dr07}
  The lisse sheaf 
$ \V_\ell$ has rank $7$  and its restriction  to 
$T_{\CC}$  is irreducible with monodromy group 
Zariski dense in $G_2(\QQ_\ell).$ Moreover, the local monodromy of 
$\V_\ell$ at $0,1,\infty$  is as follows (resp.):
$$ {\rm involutory},\quad \Unip(2)^2\oplus \Unip(3),\quad \Unip(7).$$
\end{prop}

Let $\bar{s}$ be a complex point. Note that there is an injection of the topological 
fundamental group $\pi_1(\PP^1(\CC)\setminus \{0,1,\infty\},\bar{s})$
into $\pi_1(T,\bar{s}).$ It is well known that 
$\pi_1(\PP^1(\CC)\setminus \{0,1,\infty\},\bar{s})$ is generated by 
three elements $\gamma_0,\gamma_1,\gamma_\infty,$ satisfying 
the product relation $\gamma_0\gamma_1\gamma_\infty=1.$ 
The composition 
$$ \pi_1(\PP^1(\CC)\setminus\{0,1,\infty\},\bar{s})\to \pi_1(T,\bar{s})\to G_2(\QQ_\ell)$$
defines a local system $\V_\ell^\an$ on $T(\CC).$ It follows from the 
 motivic construction of $\V_\ell$ given in \cite{DR07}, Section~3, and from 
 the comparison between \'etale 
 and singular cohomology, that there exists a natural number $N$ and 
 a local $\ZZ[\frac{1}{N}]$-system $\V^\an$ such that $\V_\ell^\an=\V^\an\otimes \QQ_\ell,$
 i.e., the images $A_0,A_1,A_\infty$ 
 of $\gamma_0,\gamma_1,\gamma_\infty$ under the monodromy
 representation of $\V^\an$ can 
 be written simultanuously over $\ZZ[\frac{1}{N}].$ It follows that there exists a natural
 number $N_1$ such that for $\ell>N_1,$ the Jordan forms of the mod-$\ell$ residual matrices 
 $\bar{A}_0,\bar{A}_1,\bar{A}_\infty$ are again of the type
involution, $\Unip(2)^2\oplus \Unip(3),$ $\Unip(7)$ (resp.).  
 
It was proven by Feit, Fong and Thompson (\cite{FeitFong}, \cite{ThompsonG2}) that 
 any triple $(g_1,g_2,g_3)$ of elements in $G_2(\FF_\ell),$ which satisfies the 
 product relation $g_1g_2g_3=1$
and whose Jordan forms are of the above type, generates the group 
$G_2(\FF_\ell)$ if $\ell>5.$ Moreover, any triple  $(g_1,g_2,g_3)$ in $\GL_7$ 
which satisfies the product relation and whose Jordan canonical forms are of type 
involution, $\Unip(2)^2\oplus \Unip(3),$ $\Unip(7)$ (resp.) is rigid, i.e., determined   
up to simultanuous conjugation by these properties (cf.~\cite{SV}). 
Therefore, we can assume that  the matrices 
$\bar{A}_0,\bar{A}_1,\bar{A}_\infty$ generate the group $G_2(\FF_\ell)$ if 
 $\ell>N_1.$ Summarizing we obtain the following result:
 
 \begin{prop}\label{dr07bis} The monodromy 
 matrices $A_0,A_1,A_\infty$ of $\V^\an_\ell$ can be written
 as elements in the group $\GL_7(\ZZ[\frac{1}{N}])$ for some $N\in \NN.$  There exists 
 a natural number $N_1>N$ such that 
for all $\ell>N_1,$ the Jordan forms of the reduction $\bar{A}_0,\bar{A}_1\bar{A}_\infty$ 
 modulo-$\ell$ of the monodromy matrices $A_0,A_1,A_\infty$ (resp.) 
  are the reduction of the respective Jordan forms. Moreover, 
  $$\langle  \bar{A}_0,\bar{A}_1,\bar{A}_\infty\rangle = G_2(\FF_\ell).$$ 
 \end{prop}

 \begin{cor}\label{dr07tris} For $\ell>N_1,$ the image of 
 the monodromy representation $\rho_{\ell}(3): \pi_1(T)\to \GL_7(\QQ_\ell)$ 
 of the third Tate twist  $\V_\ell(3)$ 
 coincides with the group $G_2(\ZZ_\ell).$ 
\end{cor}
\proof It is shown in \cite{dk} (Thm.~1) that the image of $\rho_{\ell}(3)$ is contained in 
$G_2(\QQ_\ell)$ for $\ell>2.$  
 Since the group $\pi_1(T,\bar{s})$ is compact, the image of  $\rho_{\ell}(3)$  is contained in a maximal compact subgroup
 of $G_2(\QQ_\ell).$ 
 It follows from Bruhat-Tits theory that there are $3$ (=Lie rank of $G_2$ plus $1$) distinct conjugacy
 classes of maximal compact subgroups in $G_2(\QQ_\ell).$ One is the group
 $G_2(\ZZ_\ell),$ the others are labeled by the simple roots $\Delta=\{\alpha,\beta\}$ of $G_2$
 and denoted by $G_\alpha(\ZZ_\ell)$ and $G_\beta(\ZZ_\ell)$ in \cite{GrossNebe}.  
 It is shown in \cite{GrossNebe}, Thm.~1,  that  for $\gamma=\alpha,$ or $\beta,$ the group
 $G_\gamma(\FF_\ell)$ is the semidirect product of a reduced semisimple algebraic group
 $G_\alpha^\red$
 over $\FF_\ell$ by its unipotent radical. Moreover, it is shown there that the  
 set of roots
 $$ \Phi_\gamma:=(\Delta\setminus \{\gamma\}) \cup \{-\beta_0\}$$ is a set of simple roots 
 for $G_\alpha^\red,$ where $\beta_0$ denotes the highest root of $G_2.$ 
 Therefore, we end up with $G_\alpha^\red$ to be of type $A_2$ and with 
 $G_\beta$ to be of type $A_1\cup A_1.$ Let $\ell>N_1,$ where $N_1$ is as in Prop.~\ref{dr07bis}.
 Then the reduction modulo $\ell$ of the monodromy matrices generate the group
 $G_2(\FF_\ell)$ and the residual representation of $\rho_{\ell}(3)$ is hence not 
 of type $A_2$ or $A_1\cup A_1.$ Therefore,  the image of $\rho_{\ell}(3)$ can be assumed 
 to be contained in $G_2(\ZZ_\ell).$ Since the natural surjection
 $G_2(\ZZ_\ell)\to G_2(\FF_\ell)$ has the Frattini property by \cite{Wei} (i.e., its kernel is contained 
 in the intersection of all maximal compact open subgroups),
 the image of $\rho_{\ell}(3)$ coincides with $G_2(\ZZ_\ell).$ \Endproof
 
The following result is similar to \cite{dk}, Prop.~4:
\begin{thm} \label{thmdk} Let $\ell>7$ be a prime, let $F$ be a totally real
number field and let $q,q'\not= \ell$ be odd prime numbers which split
completely in $F.$ Let $\omega,\omega'$ denote primes of $F,$
lying over $q,q'$ (resp.).  
Let $\bar{s}\in T(\bar{F})$ be a geometric point extending an $F$-rational 
point  $s\in T(F)$ such that 
$|\frac{1}{s}|_\omega<1$ and $|1-s|_{\omega'}<1.$ 
Then the following holds for the specialized Galois representation $\rho_\ell^s(3):\Gamma_F\to 
\GL_7(\QQ_\ell).$ 
\begin{enumerate}
\item 
The restriction of $\rho_\ell^s(3)$
to the inertia subgroup 
$I_\omega\leq \Gamma_{F_\omega}$ is unipotent and indecomposable,
i.e.,  of type $\Unip(7).$ If $\ell>N_1,$ where $N_1$ is the constant occurring 
in Prop.~\ref{dr07bis}, then the restriction
of $\bar{\rho}_\ell^s$ to $I_{\omega}$ is also unipotent and indecomposable. 
Moreover, semisimplification of the $\Gamma_{F_\omega}$-representation
$\rho_\ell^s(3)|_{\Gamma_{F_\omega}}$ (resp. $\bar{\rho}_\ell^s|_{G_{F_\omega}}$)
 is unramified and the eigenvalues of 
$\Frob_\omega$ are of the form $q^{-3},q^{-2},q^{-1},0,q,q^2,q^3.$
\item 
The restriction of $\rho_\ell^s(3)$ 
to the inertia subgroup 
$I_{\omega'}\leq \Gamma_{F_{\omega'}}$ is is of type $\Unip(3)\oplus \Unip(2)\oplus
\Unip(2).$ If $\ell>N_1,$ then the restriction of $\bar{\rho}_\ell^s$ 
to $I_{\omega'}$ is again of type $\Unip(3)\oplus \Unip(2)\oplus \Unip(2).$ 
\end{enumerate}
\end{thm}
\proof The proof is analogous to the proof of \cite{HSBT}, Lemma~1.15:
Let $W_q$ denote the ring of Witt vectors 
of $\bar{\FF}_q$ and $\tilde{F}$ denote its field of fractions
(the maximal unramified extension of $F_\omega$).  Let 
$t$ denote the standard parameter of $\AA^1.$
It follows from \cite{SGA1}, XIII.5.3, that there is a 
commutative diagram 
\begin{equation}\label{eqsga1} 
\begin{array}{ccc}
\pi_1(\Spec(\bar{F_\omega}((\frac{1}{t}))\,))&
 \stackrel{\sim}{\longrightarrow} &\prod_p\ZZ_p\\
\downarrow && \downarrow\\
\pi_1(W_q((\frac{1}{t})))&
\stackrel{\sim}{\longrightarrow}&
\prod_{p\not=q}\ZZ_p\\
\uparrow&&\uparrow\\
\pi_1(\Spec(\tilde{F}))&\longrightarrow &\prod_{p\not= q}\ZZ_p,
\end{array}
\end{equation}
where the left hand up-arrow is induced by $t\mapsto s,$ the right
hand downarrow is the natural projection, and the    right 
uparrow is multiplication by $\nu_{F_\omega}(s).$ 
Let $\V[\ell]$ be the lisse ${\FF}_\ell$-sheaf corresponding to the
residual representation of $\rho_\ell.$ Then 
the restriction of the sheaves 
$\V_\ell(3)$ and $\V[\ell](3)$ 
to $\Spec(W_q((\frac{1}{t})))$ correspond 
to representations of $\pi_1(\Spec(W_q((\frac{1}{t})))).$ 
It follows from Prop.~\ref{dr07} and Prop.~\ref{dr07bis} that the pullback of these representations to $\pi_1(\Spec(\bar{F_\omega}((\frac{1}{t}))))\simeq \prod_p\ZZ_p$ 
along the left downarrow sends $1$ to a unipotent matrix with minimal 
polynomial $(X-1)^7$ (perhaps enlarging $N_1$ if $\nu_{F_\omega}(t)$ is 
divided by $\ell$).  Since $\Frob_\omega$ acts on 
the inertia via the cyclotomic character and since the weight of the determinant of 
$\rho_\ell^s(3)$ is zero, the eigenvalues of 
$\Frob_\omega$ are of the given type, proving the first claim. The
second claim follows along the same arguments, using again 
Prop.~\ref{dr07} and Prop.~\ref{dr07bis}.
\Endproof

\begin{lem}\label{lemdk} Let $F$ be a number field and let  $s\in T(F).$  
\begin{enumerate} \item The system of specialized 
Galois representations $(\rho_\ell^s:\Gamma_F\to \GL_7(\QQ_\ell))_\ell$ is strictly 
compatible.
\item If $\omega$ is a prime of $F$ lying over $\ell>2,$ then the  restriction of $\rho_\ell^s$ to $\Gamma_{F_\omega}$  is de Rham
with Hodge-Tate numbers $0,1,\ldots,6.$ 
\item If $\ell$ is large enough  and if  $\omega$ is a prime of $F$ lying over $\ell,$  then the restriction of $\rho_\ell^s$ to $\Gamma_{F_\omega}$ 
is crystalline.
\end{enumerate}
\end{lem}
\proof Claim~(i) follows from \cite{Katz96}, 5.5.4. To prove the second claim we argue as follows:
By the algebraic construction of the middle convolution (cf.~\cite{dr00})
 the monodromy group of the analytification 
 $$\V^\an=\p\left[\ker\left(R^6(\pi^{\rm an}_{X})_*(\QQ) \to R^6(\pi^{\rm an}_{\coprod D_i})_*(\QQ)\right)\right]$$ 
 of $\V_\ell$  is contained in $\GL_7(\ZZ)$ and leaves 
 a bilinear form over $\ZZ$ -invariant (in the notation of Lemma~\ref{cormotg4}, where 
 $\pi^{\rm an}_X$  denotes the analytification of the map $\pi:X\to \AA^1_\QQ\setminus\{0,1\}$). 
  Hence, the analytification  $\V^\an$ of $\V$ can be seen as 
   a polarized $\ZZ$-variation of Hodge structures.
 Then the monodromy theorem (cf.~\cite{PetersSteenbrink}, Thm.~11.8) implies that the Hodge numbers of $\V^\an$ which are 
 different from zero are $h^0,h^1,\ldots,h^6.$ 
 It follows then from the  comparison isomorphism 
 that the Hodge-Tate numbers of $\rho_\ell(3)$ are 
 $-3,-2,-1,0,1,2,3.$  (The comparison isomorphism can be applied, since
the projector $\p$ in Lemma~\ref{cormotg4}~(ii) is algebraic, and hence de Rham, and
since $$\ker\left(R^6(\pi_{X})_*(\bQl) \to R^6(\pi_{\coprod D_i})_*(\bQl)\right)$$
is defined by a map between smooth projective varieties.)
The third claim
follows from the crystalline comparison isomorphism 
 since for $\ell$ large enough, the fibre $X_s$ is smooth over 
$\ZZ_\ell.$  \Endproof

\begin{thm}\label{dkkor} Suppose that $s\in T(\QQ)$ is such that there exist two 
different primes $p,q\not= 2$ such that $\nu_p(s)>0$ and $\nu_q(1-s)<0.$ Then,
for all primes  $\ell$ up to a set of density zero, the image of $\rho_\ell^s(3)$ coincides
with $G_2(\ZZ_\ell).$ 
\end{thm}

\proof It is shown in \cite{dk}, Thm.~1, that the image of $\rho_\ell^s(3)$ is Zariski dense
in $G_2(\QQ_\ell).$ By Lemma~\ref{lemdk}(i), the system $(\rho_\ell^s(3))$ is compatible. 
It follows hence from the main result in \cite{Larsen} that for all primes $\ell$ up to 
a set of density zero, the image of $\rho_\ell^s(3)$ is a maximal compact hyperspecial subgroup
of $G_2(\QQ_\ell)$ and is hence one of the groups $G_2(\ZZ_\ell),$ $G_\alpha(\ZZ_\ell)$ or 
$G_\beta(\ZZ_\ell)$ from the proof of Cor.~\ref{dr07tris}.  This implies the claim since 
the image of  $\rho_\ell^s(3)$ is contained in $G_2(\ZZ_\ell)$ by 
Cor.~\ref{dr07tris}.\Endproof

\section{Potential automorphy}\label{modularlifting}

Let us recall some definitions from \cite{BLGGT}:
Suppose that $F$ and $M$ are number fields, that $S$ is a finite set of primes of $F$ and
that $n$ is a positive integer. Then a weakly compatible system of $n$-dimensional $\ell$-adic
representations of $G_F$ defined over $M$ and unramified outside $S$ is a
family of continuous semi-simple representations
$(r_{\lambda} : G_F \to  \GL_n(M_{\lambda}))$
(where ${\lambda}$ runs over the finite places of $M$) with the following properties:

\begin{itemize}\item 
If ${\nu} \in S$ is a finite place of $F,$ then for all ${\lambda}$ not dividing the residue characteristic
of ${\nu},$ the representation $r_{\lambda}$ is unramified at ${\nu}$ and the characteristic
polynomial of $r_{\lambda}(\Frob_{\nu})$ lies in $M[X]$ and is independent of ${\lambda}.$
\item Each representation $r_{\lambda}$ is de Rham at all places above the residue characteristic
of ${\lambda},$ and in fact crystalline at any place ${\nu} \notin S$ which divides the
residue characteristic of ${\lambda}.$
\item For each embedding $\tau : F \to \bar{M},$  the $\tau$-Hodge-Tate numbers of $r_{\lambda}$ are
independent of ${\lambda}.$
\end{itemize}

Let us recall the potential automorphy criterion of \cite{BLGGT} for 
weakly compatible systems (\cite{BLGGT}, Thm.~A and Thm.~5.3):

\begin{thm}\label{AL} Let $F$ be a totally real number field.
Let $(r_\lambda)$ be a weakly compatible system of $n$-dimensional $\ell$-adic 
representations of $G_F$ defined over $M$ and unramified outside $S,$ where for simplicity
we assume that $M$ contains the image of each embedding $F \to M.$ Suppose that
$(r_{\lambda})$ satisfies the following properties.
\begin{enumerate}
\item (Irreducibility) There exists a set of primes $L$ of density $1$ such that 
$r_{\lambda}$ ($\lambda \in L$)  is irreducible.
\item (Hodge-Tate Regularity) For each embedding $\tau : F \to M$ the representation $r_{\lambda}$ has $n$
distinct $\tau$-Hodge-Tate numbers.
\item (Odd essential self-duality) Either each $r_{\lambda}$ factors
through a map to ${\rm GSp}_n(M_{\lambda})$ with a totally odd multiplier character; or
each $r_{\lambda}$ factors through a map to $GO_n(M_{\lambda})$ with a totally even multiplier
character. Moreover in either case the multiplier characters form a weakly
compatible system.
\end{enumerate}
Then there is a finite, Galois, totally real extension over which all the $r_{\lambda}$Õs become
automorphic. \end{thm}

\begin{thm}\label{thmhaupt} Let $s\in T(\QQ)$ is such that there exist two 
different primes $p,q\not= 2$ such that $\nu_p(s)>0$ and $\nu_q(1-s)<0.$
Then  the  Galois 
representations $\rho_\ell^s(3)$ are potentially automorphic for all 
$\ell.$ 
\end{thm}

\proof By Lemma~\ref{lemdk} (i), we known that the system of Galois representations
$(\rho_\ell^s(3))$ is compatible. By Theorem~\ref{dkkor}, the system $(\rho_\ell^s(3))$
is irreducible for a set $L$ of primes having density $1.$ 
Moreover, since each $\rho_\ell^s$ factors through the simple group $\G_2\leq \SO_7$ and so, 
$\rho_\ell^s(3)$ has a totally even multiplier character. By Lemma~\ref{lemdk}, we can
assume (by possibly shrinking $L$) that $\rho_\lambda$ is crystalline with Hodge-Tate numbers 
$-3,-2,-1,0,1,2,3.$ Hence we can apply Thm.~\ref{AL}, the claim follows.\Endproof

\bibliographystyle{plain} \bibliography{/users/michaeldettweiler/michael/Biblio/p}

\end{document}